\documentclass[12pt]{amsart}
\usepackage{amsmath,amsfonts,mathrsfs,amssymb,cite}
\usepackage{bm}

\newtheorem{thm}{Theorem}[section]
\newtheorem{theorem}[thm]{Theorem}
\newtheorem{corollary}[thm]{Corollary}
\newtheorem{lemma}[thm]{Lemma}
\newtheorem{proposition}[thm]{Proposition}

\theoremstyle{definition}

\theoremstyle{remark}

\newenvironment{theorem*}[1]{\smallskip\noindent{\bf #1.}\it}{\medskip}

\numberwithin{equation}{section}

\newcommand\tr{\operatorname{tr}}
\newcommand{\slim}{\operatornamewithlimits{s-lim}}
\newcommand{\re}{\mathrm{e}}
\newcommand{\ri}{\mathrm{i}}
\newcommand{\diag}{\operatorname{diag}}
\newcommand{\dom}{\operatorname{dom}}

\newcommand{\bC}{{\mathbb C}}
\newcommand{\bL}{{\mathbb L}}
\newcommand{\bN}{{\mathbb N}}
\newcommand{\bH}{{\mathbb H}}
\newcommand{\bR}{{\mathbb R}}
\newcommand{\bZ}{{\mathbb Z}}

\newcommand\cH{{\mathcal H}}

\newcommand\fI{{\mathfrak I}}
\newcommand\fH{{\mathfrak H}}
\newcommand\fS{{\mathfrak S}}
\newcommand\fP{{\mathfrak P}}

\newcommand{\sA}{{\mathscr A}}
\newcommand{\sB}{{\mathscr B}}

\newcommand{\sD}{{\mathscr D}}
\newcommand{\sE}{{\mathscr E}}
\newcommand{\sH}{{\mathscr H}}
\newcommand{\sI}{{\mathscr I}}
\newcommand{\sF}{{\mathscr F}}
\newcommand{\sK}{{\mathscr K}}
\newcommand{\sL}{{\mathscr L}}

\newcommand{\sN}{{\mathscr N}}
\newcommand{\sR}{{\mathscr R}}
\newcommand{\sS}{{\mathscr S}}
\newcommand{\sT}{{\mathscr T}}
\newcommand{\sX}{{\mathscr X}}
\newcommand{\sY}{{\mathscr Y}}

\newcommand\al{\alpha}
\newcommand\be{\beta}
\newcommand\ga{\gamma}
\newcommand\la{\lambda}
\newcommand\eps{\varepsilon}
\newcommand\si{\sigma}

\newcommand{\bs}{{\mathbf s}}

\newcommand{\bu}{{\mathbf u}}
\newcommand{\bx}{{\mathbf x}}
\newcommand{\by}{{\mathbf y}}

\newcommand\bal{\bm{\alpha}}
\newcommand\bga{\bm{\gamma}}
\newcommand\bde{\bm{\delta}}

\newcommand\bla{\bm{\lambda}}
\newcommand\bmu{\bm{\mu}}
\newcommand\bnu{\bm{\nu}}

\begin{document}

\title[Uniform stability]{Analyticity and uniform stability in
     the inverse spectral problem for Dirac operators}%
\author{Rostyslav O. Hryniv}%
\address{Institute for Applied Problems of Mechanics and Mathematics, 3b Naukova st., 79601 Lviv, Ukraine
\and Institute of Mathematics, the University of Rzesz\'{o}w, 16\,A Rejtana al., 35-959 Rzesz\'{o}w, Poland}
\email{rhryniv@iapmm.lviv.ua}

\thanks{}%

\date{12 February 2011}%

\begin{abstract}
We prove that the inverse spectral mapping reconstructing the square
integrable potentials on~$[0,1]$ of Dirac operators in the AKNS form
from their spectral data (two spectra or one spectrum and the
corresponding norming constants) is analytic and
uniformly stable in a certain sense.
\end{abstract}

\maketitle



\section{Introduction}

The main goal of this paper is to establish analyticity and uniform continuity of solutions in the inverse spectral problems for a certain class of Dirac operators on the interval~$[0,1]$. Namely, we consider one-dimensional Dirac operators in the AKNS normal
form~\cite{AKNS} generated by the differential expressions
\[
    \ell_Q:= \sigma_2\frac1{\ri}\frac{d}{dx} +Q(x),
\]
with potentials
\begin{equation}\label{eq:intr.Q}
    Q(x) = q_1(x) \sigma_1 + q_3(x) \sigma_3.
\end{equation}
Here
\[
    \sigma_1 = \begin{pmatrix} 0 &1 \\ 1 & 0
               \end{pmatrix}, \qquad
    \sigma_2 = \begin{pmatrix} 0 &-\ri \\ \ri & 0
               \end{pmatrix}, \qquad
    \sigma_3 = \begin{pmatrix} 1& 0  \\ 0 & -1
               \end{pmatrix}
\]
are the Pauli matrices and $q_1$ and $q_2$ are real-valued integrable
functions that incorporate various physical characteristics of the
particle considered (its mass $m$, the scalar potential, electrostatic
potential, and anomalous magnetic moment, see~\cite[Ch.~4]{Th}). We
assume that the potential is of bounded support and thus consider the
Dirac operator on a bounded interval (which without loss of generality
might be taken $[0,1]$).

The Dirac equation has widely been used in various areas of physics and mathematics starting from 1929, when P.~Dirac suggested it to model the evolution of spin-$\tfrac12$ particles in
the relativistic quantum mechanics~\cite{Th}. In 1973 Ablowitz, Kaup,
Newell, and Segur~\cite{AKNS} observed that the Dirac equation
is related to a nonlinear wave equation (the ``modified
Korteweg--de Vries equation", a member of the AKNS--ZS hierarchy,
see~\cite{AKNS2,ZS}) in the same manner as the Schr\"odinger
equation is related to the KdV equations, and this discovery stimulated the  increasing interest in direct and inverse problems for Dirac operators in both physical and mathematical literature. The topic was brought to life somewhat earlier; namely, in 1966, Gasymov and Levitan solved the inverse problems for Dirac operators on $\bR_+$ by the spectral function~\cite{GasLsp} and by
the scattering phase~\cite{GasLscat}. Their investigations were
continued and further developed in many directions. The reference list on the inverse spectral and scattering theory for Dirac operators has become very vast by now, and we only mention papers~\cite{AHMdir,CG,HJKS,Mal99} that contain extended bibliography on the topic. The books
by Levitan and Sargsjan~\cite{LS} and by Thaller~\cite{Th} may
serve as a good introduction into the (respectively mathematical
and physical part of the) theory of Dirac operators.

In what follows, we fix $\theta_0,\theta_1\in[0,\pi)$ and define the Dirac operator
$\sD(\theta_0,\theta_1,Q)$ in the Hilbert space $\bH:=L_2(0,1)\times
L_2(0,1)$ as the restriction of $\ell_Q$ onto the domain
\begin{multline*}
    \dom \sD(\theta_0,\theta_1,Q) = \bigl\{\bu=(u_1,u_2)^{\mathrm{t}}\in \bH \mid
        u_1,u_2 \in \mathrm{AC}[0,1], \\
        u_1(j)\cos\theta_j = u_2(j)\sin\theta_j, j=0,1
            \bigr\}.
\end{multline*}
It is well known that so defined operator $\sD(\theta_0,\theta_1,Q)$ is
self-adjoint and has a simple discrete spectrum tending to $\pm\infty$.
In general, the spectrum of the operator~$\sD(\theta_0,\theta_1,Q)$ does not determine the potential~$Q$.
Indeed, there are only several analogues of the Ambartsumyan result for the Sturm--Liouville operators on a finite interval~\cite{Am}, where the unperturbed Dirac operator (i.e., with $Q=0$) is the only one having the specific spectrum; see~\cite{Ho,Ki,Yang}. To reconstruct the operator~$\sD(\theta_0,\theta_1,Q)$ in a generic situation some additional information is needed. For instance, the classical result of the inverse theory for the Dirac operators states that if we also know the spectrum of another operator $\sD(\theta_0,\theta_1',Q)$, with $\theta_1'$ such that $\sin(\theta_1'-\theta_1)\ne0$, then two such spectra determine uniquely the potential $Q$ in the AKNS form~\eqref{eq:intr.Q}.

For the sake of definiteness (and in order to avoid inessential
technicalities), we choose the boundary conditions corresponding to
$\theta_0=\theta_1=0$ and  $\theta_1'=\pi/2$ and denote the operators
$\sD(0,0,Q)$ and $\sD(0,\pi/2,Q)$ by $\sD_1(Q)$ and $\sD_2(Q)$
respectively. If $Q\in L_2(0,1)\otimes\bC^4$, then the eigenvalues $(\la_n)_{n\in
\bZ}$ and $(\mu_n)_{n\in \bZ}$ of $\sD_1(Q)$ and $\sD_2(Q)$
respectively can be enumerated so that they satisfy the interlacing
condition
\begin{equation}\label{eq:intr.interlace}
\la_n < \mu_n < \la_{n+1}, \qquad n\in\bZ,
\end{equation}
and the asymptotics
\begin{equation}\label{eq:intr.as-la}
 \begin{aligned}
    \la_n  &= \pi n +  \la'_n,\\
    \mu_n  &= \pi (n+\tfrac12) +  \mu'_n,
 \end{aligned}
\end{equation}
where the remainders $\la'_n$ and $\mu'_n$ form $\ell_2$-sequences; see~\cite{Amour,LS}.
Conversely, it was shown in~\cite{AHMdir} that any two sequence of real numbers that interlace and obey the asymptotics~\eqref{eq:intr.as-la} are eigenvalues of the operators $\sD_1(Q)$ and $\sD_2(Q)$, for some real-valued potential~$Q$ in the AKNS form with $L_2$-entries.

We thus have a well-defined mapping
\begin{equation}\label{eq:intr.spmap}
    \bigl((\la_n),(\mu_n)\bigr) \mapsto Q
\end{equation}
from the spectral data~$\bigl((\la_n),(\mu_n)\bigr)$ of Dirac operators $\sD_1(Q)$ and $\sD_2(Q)$ to the corresponding potential $Q$ in the AKNS form. A natural question arises, what properties this inverse spectral mapping has.

Unfortunately, we are not aware of any work where this topic has been explicitly studied. However, there have been many papers discussing the analogous question for Sturm--Liouville operators with integrable potentials; see~\cite{Rya3,McL,An,CM} and the references cited in~\cite{Hstab}. (We observe in passing that the case of Sturm--Liouville operators in impedance form considered in~\cite{McL,An,CM,Hstab-imp} is in a certain sense equivalent to the class of reduced Dirac operators with $q_3=0$; cf.~\cite{AHMdir}.) Adapting as necessary the corresponding proofs for the Sturm--Liouville operators, one can show that the inverse spectral problem is locally stable, i.e., that the potential~$Q$ depends continuously on the spectral data $\bigl((\la_n),(\mu_n)\bigr)$. Here the metric on the set of the spectral data is induced by that of $\ell_2\times\ell_2$ on the elements $\bigl((\la'_n),(\mu'_n)\bigr)$ with the $\la_n'$ and $\mu_n'$ defined by~\eqref{eq:intr.as-la}. More exactly (cf.~\cite{McL,An} for the Sturm--Liouville case), this local stability means that, for a fixed $M>0$, there are positive~$\eps$ and $L$ such that if $\|Q_1\|_{\bH}\le M$, $\|Q_2\|_{\bH}\le M$, and their spectral data $\bnu_1:=\bigl((\la_{1,n}),(\mu_{1,n})\bigr)$ and
$\bnu_2:=\bigl((\la_{2,n}),(\mu_{2,n})\bigr)$ satisfy
$\|\bnu_1-\bnu_2\|\le\eps$, then
\begin{equation}\label{eq:intr.cont}
    \|Q_1-Q_2\|_{\bH} \le L \|\bnu_1-\bnu_2\|.
\end{equation}

Such a result cannot be considered satisfactory e.g.\ for the purpose
of numerical reconstruction, as it refers to the norm of the
potential~$Q$ to be recovered and thus specifies neither the allowed
noise level~$\eps$ nor the Lipschitz constant~$L$. Therefore it is desirable to have a uniform stability that asserts~\eqref{eq:intr.cont} whenever the spectral data $\bnu_1$ and $\bnu_2$ run through a bounded set~$\sN$ and with $L$ only depending on~$\sN$.

Recently, such a uniform stability in the inverse spectral problem for
Sturm--Liouville operators on~$[0,1]$ was established by Shkalikov and Savchuk~\cite{SSstab}. They considered operators with real-valued potentials from the Sobolev spaces~$W_2^{s}(0,1)$ with $s>-1$; for negative $s$, such potentials are distributions. Their approach was based on a careful study of the nonlinear mapping of~\eqref{eq:intr.spmap} and used extensively the implicit function theorem. In our work~\cite{Hstab} we established analyticity and global stability of the inverse spectral mapping for $s\in[-1,0]$ using a different approach that generalizes a classical method due to Gelfand and Levitan~\cite{GL} and Marchenko~\cite{Ma} and has been successfully applied to reconstruction of Sturm--Liouville operators with singular potentials~\cite{HMinv,HMscale}. In~\cite{Hstab-imp}, similar results were established for Sturm--Liouville operators in impedance form.

The main aim of this paper (Theorem~\ref{thm:pre.main}) is to prove analyticity and uniform continuity of the inverse spectral mapping~$\bigl((\la_n),(\mu_n)\bigr) \mapsto Q$ for the class of the Dirac operators under consideration. In fact, we first establish these properties for the inverse problem of reconstructing $Q$ from the eigenvalues of~$(\la_n)$ and the corresponding norming constants~$(\alpha_n)$ introduced below (Theorem~\ref{thm:pre.norm}), and then show that these norming constants depend analytically an Lipschitz continuously on two spectra, $(\la_n)$ and $(\mu_n)$. We use the approach to the inverse spectral problem for Dirac equation based on the Krein equation~\cite{AHMdir} and further develop the methods of~\cite{Hstab}.

We mention that the methods of~\cite{AHMdir} could be used to treat the Dirac operators in AKNS form with real-valued potentials belonging to $L_p(0,1)\otimes \bC^4$ with $p\in[1,\infty)$. However, apart from some technicalities caused by more complicated properties of the Fourier transform in $L_p(0,1)$ for $p\ne2$ the approach would remain the same and we decided to sacrifice the generality to simplicity of presentation.

The paper is organised as follows. In the next section we introduce all the related objects, formulate the main results, and recall the approach to solution of the inverse spectral problem for Dirac operators in AKNS form based on the Gelfand--Levitan--Marchenko and Krein equations. In Section~\ref{sec:stab-norm}, we prove Theorem~\ref{thm:pre.norm} and in Section~\ref{sec:two} study dependence of the spectrum $(\mu_n)$ of the operator~$\sD_2(Q)$ on $(\la_n)$ and $(\al_n)$ and then derive Theorem~\ref{thm:pre.main}. Finally, three appendices discuss properties of special analytic mappings between Banach spaces and recall some facts on the special Banach algebra that are essentially used in the paper.

\emph{Notations.} Throughout the paper, we shall denote by $\bL_2(0,1)$ the set of all $2\times2$ potentials in the AKNS normal form; clearly, this space is unitarily equivalent to $L_2(0,1) \times L_2(0,1)$. Bold letters will usually denote vectors and functions with values in~$\bH$, and $\|f\|$ and $\|\mathbf{f}\|$ will stand for the $L_2(0,1)$-norm of a scalar function~$f$ and $\bH$-norm of the vector-valued function~$\mathbf{f}$, respectively.


\section{Preliminaries and main results}\label{sec:pre}

In this section we state the main results and recall the method of
solution of the inverse spectral problem based on the
Gelfand--Levitan--Marchenko~\cite{LS} and Krein~\cite{AHMdir} equations.

\subsection{Spectral data}

We denote by $\la_n$ and $\mu_n$, $n\in\bZ$, the eigenvalues of the
operators $\sD_1(Q)$ and $\sD_2(Q)$ respectively and recall that these
eigenvalues interlace, i.e., $\la_n<\mu_n< \la_{n+1}$ for all
$n\in\bZ$, and satisfy the relations
\[
    \la_n = \pi n + \rho_{2n},
    \qquad
    \mu_n = \pi (n+\tfrac12) + \rho_{2n+1}
\]
with some $\ell_2$-sequence $(\rho_n)$.

The equation $\ell_Q \bu = \la \bu$, $\bu=(u_1,u_2)^\mathrm{t}$, subject to
the initial conditions $u_1(0)=0$, $u_2(0)=1$ has the solution
\begin{equation}\label{eq:pre.s}
    \bs(x,\la) =
        \bs_0(x,\la)
        + \int_0^x K(x,t) \bs_0(t,\la)\,ds,
\end{equation}
where $\bs(x,\la):=\bigl(s_1(x,\la),s_2(x,\la)\bigr)^\mathrm{t}$, $\bs_0(x,\la):=(\sin\la x, \cos\la x)^\mathrm{t}$ and $K=:(k_{jl})_{j,l=1}^2$ is the kernel of the transformation
operator.
Further, $\bs(\cdot,\la_n)$ is an eigenfunction corresponding to the
eigenvalue~$\la_n$ of the operator~$\sD_1(Q)$, and we call the number
$\al_n := \|\bs(\cdot,\la_n)\|^{-2}$ the \emph{norming constant} for
this eigenvalue. It is known that
\[
    \al_n = 1 + \alpha'_n,
\]
where the sequence $(\alpha'_n)_{n\in\bZ}$ belongs to $\ell_2:=\ell_2(\bZ)$.
Moreover, the norming constants~$\al_n$ can be determined from the
spectra of the operators~$\sD_1(Q)$ and $\sD_2(Q)$ as follows.

We set $S(\la):=s_1(1,\la)$ and $C(\la):=s_2(1,\la)$; due
to~\eqref{eq:pre.s} these are entire functions of exponential type~$1$
with zeros $\la_n$ and $\mu_n$ respectively. The Hadamard canonical
products of~$S$ and $C$ are
\begin{equation}\label{eq:pre.PhiPsi}
  S(z) = (z-\la_0) \mathrm{V.p.}\hspace{-3pt}
        \prod\limits_{n=-\infty}^{\infty}
        \hspace{-7pt}{\vphantom{\prod}}^\prime
        \hspace{5pt}\frac{\la_n-z}{\pi n},
        \qquad
  C(z) = \mathrm{V.p.}\hspace{-3pt}
        \prod\limits_{n=-\infty}^{\infty}
        \frac{\mu_n -z}{\pi(n+\tfrac12)}
\end{equation}
(the prime denoting that the factor corresponding to $n=0$ is omitted),
so that $S$ and $C$ are uniquely determined by their zeros. Then
we have~\cite{AHMdir}
\begin{equation}\label{eq:pre.al}
        \al_n = \frac{1}{\dot{S}(\la_n)C(\la_n)},
\end{equation}
where the dot denotes the derivative in~$z$.


\subsection{Main results}
We introduce the set~$\sN$ of data
$\bigl((\la_n)_{n\in\bZ},(\mu_n)_{n\in\bZ}\bigr)$ with the following
properties:
\begin{itemize}
\item the sequences $(\la_n)$ and $(\mu_n)$ strictly interlace,
    i.e., $\la_n<\mu_n<\la_{n+1}$ for all $n\in\bZ$;
\item the sequence $(\rho_k)_{k\in\bZ}$, with
    $\rho_{2n}:= \la_n - \pi n$ and
    $\rho_{2n+1}:=\mu_n- \pi (n+\tfrac12)$, belongs to $\ell_2$.
\end{itemize}
In this way every element $\bnu:=\bigl((\la_n),(\mu_n)\bigr)$ of $\sN$ is
identified with a sequence $(\rho_n)$ in~$\ell_2$ thus inducing a
metric on~$\sN$.

According to~\cite{AHMdir}, every element of~$\sN$ gives the eigenvalue
sequences of the operators $\sD_1(Q)$ and $\sD_2(Q)$ corresponding to
a unique real-valued AKNS potential~$Q\in \bL_2(0,1)$ and,
conversely, for every real-valued $Q\in \bL_2(0,1)$ of the
form~\eqref{eq:intr.Q} the spectra of the corresponding Dirac operators
form an element of~$\sN$. When the AKNS potential~$Q$ varies over a
bounded subset of~$\bL_2(0,1)$, then the corresponding
spectral data $\bigl((\la_n),(\mu_n)\bigr)$ remain in a bounded subset
of $\sN$. Moreover, the Pr\"ufer angle technique yields then a positive
$h$ such that all the corresponding spectral
data~$\bigl((\la_n),(\mu_n)\bigr)$ are $h$-se\-pa\-rated, i.e., that
$\mu_{n+1}-\la_n\ge h$ and $\la_n-\mu_n\ge h$ for every $n\in\bN$.
Summarizing, we conclude that the uniform stability of the inverse
spectral problem we would like to establish is only possible on the
convex closed sets $\sN(h,r)$ of spectral data consisting of all
elements of~$\sN$ that are $h$-separated and satisfy $\|(\rho_n)\|_{\ell_2}\le r$.

In these notations, the first main results of the paper reads
as follows.

\begin{theorem}\label{thm:pre.main}
For every $h\in(0,\pi/2)$ and $r>0$, the inverse spectral mapping
\[
    \sN(h,r) \ni \bnu \mapsto Q \in \bL_2(0,1)
\]
is analytic and Lipschitz continuous.
\end{theorem}

In fact, as in~\cite{Hstab}, we prove first the analyticity and Lipschitz continuity of the
inverse spectral problem of reconstructing $Q$ from the Dirichlet
spectrum $(\la_n)$ and the norming constants $(\al_n)$ (see
Theorem~\ref{thm:pre.norm} below), and then derive
Theorem~\ref{thm:pre.main} by showing that the norming constants depend
analytically and Lipschitz continuously on the two spectra.

More exactly, we denote by $\sL$ the family of strictly increasing
sequences $\bla:=(\la_n)$ such that $\rho_{2n}:=\la_n-\pi n$ forms an
element of $\ell_2$ and pull back the topology on~$\sL$ from that of $\ell_2$ by identifying
such $\bla$ with $(\rho_{2n})\in\ell_2$. For $h\in(0,\pi)$ and
$r>0$, we denote by~$\sL(h,r)$ the closed convex subset of $\sL$
consisting of sequences $(\la_n)_{n\in\bZ}$ with $\la_{n+1}-\la_n \ge
h$, and such that $\|(\rho_{2n})\|_{\ell_2}\le r$. Next, we write $\sA$ for
the set of sequences $\bal:=(\al_n)_{n\in\bN}$ of positive numbers
$(\al_n)$ such that the sequence~$(\be_n)$ with $\be_n:=\al_n-1$
belongs to~$\ell_2$. This induces the topology of~$\ell_2$ on $\sA$; we
further consider closed subsets $\sA(h,r)$ of $\sA$ consisting of all
$(\al_n)$ such that~$\al_n \ge h$ for all $n\in\bN$ and $\|(\be_n)\|_{\ell_2}\le r$.

It is known that, given an element $(\bla,\bal)\in\sL \times \sA$,
there is a unique $Q\in\bL_2(0,1)$ in the AKNS normal form such that $\bla$ is the
sequence of eigenvalues and $\bal$ the sequence of norming constants
for the Dirac operator~$\sD_1(Q)$. Some further properties of the
induced mapping are described in the following theorem.

\begin{theorem}\label{thm:pre.norm}
For every $h\in(0,\pi)$ and every positive $h'$, $r$, and $r'$, the inverse spectral
mapping
\[
    \sL(h,r) \times \sA(h',r')\ni (\bla,\bal)
        \mapsto Q \in \bL_2(0,1)
\]
is analytic and Lipschitz continuous.
\end{theorem}


\subsection{Solution of the inverse spectral problem via the GLM and Krein equations}

The Gelfand--Levitan--Marchenko equation relates the spectral data
for the operator $\sD_1(Q)$ (i.e., its eigenvalues and norming
constants) with the transformation operator $\sI+\sK$. To derive it,
we start with the resolution of identity for $\sD_1(Q)$,
\[
    \sI = \slim_{N\to\infty} \sum_{n=-N}^N \al_n (\cdot, \bs_n)_{\bH}\bs_n,
\]
where $\slim$ stands for the limit in the strong operator topology
and $\bs_n(x):= \bs(x,\la_n)$. Recalling that $\bs_n =
(\sI+\sK)\bs_{0,n}$ with $\bs_{0,n}(x) = \bs_0(x,\la_n)$, we get
\begin{equation}\label{eq:inv.I}
    \sI = (\sI+\sK)
    \Bigl[\slim_{N\to\infty} \sum_{n=1}^N
                \al_n (\cdot,\bs_{0,n})_{\bH}\bs_{0,n}\Bigr]
        (\sI+\sK^*).
\end{equation}
The operator in the square brackets has the form $\sI+\sF$, where
$\sF$ is an integral operator of Hilbert--Schmidt class with kernel
\begin{equation}\label{eq:inv.F}
    F (x,t) := \tfrac12 \bigl[ H(\tfrac{x-t}2) -
                   H(\tfrac{x+t}2)\si_3\bigr],
\end{equation}
and
\begin{equation}\label{eq:inv.H}
    H(s):= \mathrm{V.p.}\sum_{n=-\infty}^\infty
        \Bigl(\al_n \re^{2 \la_n\ri s\si_2} - \re^{2\pi n\ri s\si_2}\Bigr).
\end{equation}
Applying $(\sI+\sK^*)^{-1}$ to both sides of~\eqref{eq:inv.I} and
rewriting the resulting relation in terms of the kernels $K$ and
$F$, we get the Gelfand--Levitan--Marchenko (GLM) equation
\begin{equation}\label{eq:inv.GLM}
    K(x,t) + F(x,t)
        + \int_0^x K(x,s) F(s,t)\,ds = 0,
    \qquad x\ge y.
\end{equation}

If the potential $Q$ is continuous, then the algorithm for
reconstructing $Q$ from the spectral data $\bigl((\la_n),(\mu_n)\bigr)$
proceeds as follows, cf.~\cite{LS}. We first calculate numbers $\al_n$
via~\eqref{eq:pre.al}, then construct a matrix-function~$H$
of~\eqref{eq:inv.H}, form the kernel $F$ of~\eqref{eq:inv.F}, solve the
GLM equation~\eqref{eq:inv.GLM}, and finally set
\begin{equation}\label{eq:inv.QKB}
    Q(x) = K(x,x)\frac1{\ri}{\si_2} - \frac1{\ri}{\si_2} K(x,x).
\end{equation}

If the potential~$Q$ is only assumed to be in $\bL_2(0,1)$,
then the kernel $K$ has the property that $K(x,\cdot)$ and $K(\cdot,x)$
are square-integrable matrix-functions depending continuously in
$L_2(0,1)\otimes\bC^4$ on $x\in [0,1]$. In particular, $K$ may not have
well defined restriction $K(x,x)$ on the diagonal and thus the above-described method cannot be used.

To overcome this obstacle, we consider, along with the GLM equation, the
Krein equation
\begin{equation}\label{eq:inv.krein}
    R(x,t) + H(x-t) + \int_0^x R(x,s) H(s-t)\,ds =0,
        \qquad 0\le t\le x \le 1.
\end{equation}
It is easily seen that if $R$ is a solution to this equation, then
\[
    K(x,t):= \tfrac12 \bigl[R(x,\tfrac{x+t}{2})
            -  R(x,\tfrac{x-t}{2})\si_3\bigr]
\]
solves the GLM equation, and in this sense the Krein equation is more
general. The kernel $R$ generates another transformation operator in the sense that
\begin{equation}\label{eq:pre.sR}
   \bs(x,\la) = \bs_0(x,\la)
        + \int_0^x R(x,x-t) \bs_0(x-2t,\la)\,ds.
\end{equation}
Moreover, as soon as the related convolution operator
\begin{equation}\label{eq:pre.conv}
    \sH \bu (x) := \int_0^1 H(x-t)\bu(t)\,dt
\end{equation}
satisfies the inequality~$\sI + \sH >0$ in $L_2(0,1)\times\bC^4$, the
Krein equation can be shown to possess a unique solution $R$ and this
solution has the property that $R(\cdot,t)$ belongs to
$L_2(0,1)\otimes\bC^4$ and depends continuously therein on $t\in[0,1]$.
Now, rewriting the relation~\eqref{eq:inv.QKB} in terms of $R$, we find
that
\begin{equation}\label{eq:inv.QR}
    Q(x)  =  - R(x,0)\si_1.
\end{equation}
This formula defines $Q$ as an $L_2(0,1)\otimes\bC^4$-valued function
and will be the basis of our reconstruction algorithm.


\section{Stability of the inverse spectral problem: norming
constants}\label{sec:stab-norm}

In this section, we prove Theorem~\ref{thm:pre.norm} on uniform
stability of the inverse spectral problem of reconstructing the
potential~$Q$ of a Dirac operator $\sD_1(Q)$ from its spectrum and norming
constants.

We shall study the correspondence between the data
$(\bla,\bal)\in\sL(h,r)\times \sA(h',r')$ and the potentials $Q$ in the
AKNS form of the Dirac operator $\sD_1(Q)$ through the chain
\[
    (\bla,\bal) \mapsto H
         \mapsto R \mapsto Q,
\]
in which $H$ is the matrix-valued function of~\eqref{eq:inv.H}, $R$ is
the matrix-function solving the Krein equation~\eqref{eq:inv.krein},
and, finally, $Q$ is given by~\eqref{eq:inv.QR}.



\begin{lemma}\label{lem:phicont}
The mapping
\[
    \sL(h,r)\times \sA(h',r') \ni (\bla,\bal) \mapsto
     H\in L_2(-1,1)\otimes\bC^4
\]
is analytic and Lipschitz continuous.
\end{lemma}

\begin{proof}
Since the matrix $B$ has simple eigenvalues $\pm\ri$, it is unitarily
equivalent to the matrix $\diag\{\ri,-\ri\}$, and thus $H$ is unitarily
equivalent to the matrix-function $\diag\{h(s), h(-s)\}$, where
\begin{equation}\label{eq:stn.h}
    h(s):= \mathrm{V.p.}\sum_{n=-\infty}^\infty
        (\al_n \re^{2\la_n \ri s} - \re^{2\pi n\ri s}).
\end{equation}
We thus need to prove that the mappings
\[
    \sL(h,r)\times \sA(h',r') \ni (\bla,\bal) \mapsto h(\pm s) \in
    L_2(-1,1)
\]
are analytic and Lipschitz continuous. It clearly suffices to consider
only the case of~$h(s)$; also, since the $\al_n$ and $\la_n$ are real, we get $h(-s) = \overline{h(s)}$, and only the restriction of~$h$ onto $(0,1)$  need to be studied.

We have
 \(
     h = h_{\bla} + h_{\bla,\bal},
 \)
where
\[
    h_{\bla}(s) := \mathrm{V.p.}\sum_{n=-\infty}^\infty
        [\re^{2\rho_{2n} \ri s} - 1]\re^{2\pi n\ri s},
        \qquad
    h_{\bla,\bal}(s):= \mathrm{V.p.}\sum_{n=-\infty}^\infty
        \beta_n \re^{2\la_n \ri s}.
\]
We recall that the numbers $\beta_n:=\alpha_n-1$ and $\rho_{2n}=\la_n
- \pi n$ form sequences in~$\ell_2$ that induce the topology
of $\sL$ and $\sA$. The series
\[
    \mathrm{V.p.}\sum_{n\in\bZ} \rho_{2n}\re^{2\pi n\ri s}
    \qquad \text{and} \qquad
    \mathrm{V.p.}\sum_{n\in\bZ} \beta_n \re^{2\pi n \ri s}
\]
converge in~$L_2(0,1)$, and we denote by $f_{\bla}$ and $g_{\bal}$ respectively the
corresponding sums. Since the sequence of functions
 $\bigl(\re^{2\pi n\ri s}\bigr)_{n\in\bZ}$
forms an orthonormal basis of $L_2(0,1)$, the
mappings~$(\rho_{2n})\mapsto f_{\bla}$ and $(\beta_n)\mapsto g_{\bal}$ are unitary
isomorphisms between $\ell_2(\bZ)$ and $L_2(0,1)$ given by the inverse Fourier transform. Therefore
we have $h_{\bla} = g(f_{\bla})$ and $h_{\bla,\bal} = \Phi(f_{\bla},g_{\bal})$ with the mappings $g(f)$ and $\Phi(f,g)$ introduced in Appendix~\ref{app:aux}. The claim now follows from the properties of $g(f)$ and $\Phi(f,g)$ established in Lemmata~\ref{lem:B.gf} and \ref{lem:B.Phi} respectively.
\end{proof}

Solubility of the Krein equation crucially relies on the following property of the convolution operator~$\sH$.

\begin{lemma}\label{lem:5.Iholds}
For every $h\in(0,\pi)$, $h'\in(0,1)$, and positive $r$ and $r'$, there exists $\eps>0$ with the following property: if
$(\bla,\bal)$ is an arbitrary element of~$\sL(h,r)\times\sA(h_1,r_1)$ and $H$ is the function
of~\eqref{eq:inv.H}, then for the corresponding convolution
operator~$\sH$ we have $\sI+\sH\ge \eps\sI$.
\end{lemma}

\begin{proof}
Denote by $\sH_\pm$ the convolution operators in $L_2(0,1)$ constructed as in~\eqref{eq:pre.conv} but for the scalar functions $h(\pm s)$ of~\eqref{eq:stn.h}. Transforming~$H$ to the diagonal form as in the proof of the previous lemma, we see that it suffices to prove that there is
$\eps>0$ independent of $(\bla,\bal)$ such that $I + \sH_\pm\ge \eps I$. We shall only treat the case of the operator~$I+\sH_+$, as the other case is completely analogous.

Since the system~$\bigl(\re^{2\pi n\ri s}\bigr)_{n\in\bZ}$ is an orthonormal
basis of $L_2(0,1)$, we find that
\begin{align*}
     ((I+\sH_+) f,f)
        &= (f,f) + \lim_{k\to\infty}
            \sum_{n=-k}^k \bigl[\al_n |(f, \re^{2\la_n\ri s})|^2 -
            |(f, \re^{2\pi n\ri s})|^2\bigr]\\
        &= \mathrm{V.p.}\sum_{n=-\infty}^\infty \al_n |(f, \re^{2\la_n\ri s})|^2.
\end{align*}
The inclusion $\bal\in\sA(h',r')$ implies that $\al_n \ge h'$ for all $n\in\bZ$. It follows from the results of~\cite[Ch.~VI]{GK1}, \cite[Ch.~4]{Young} that the system~$\sE_{\bla}:=\bigl(\re^{2\la_n\ri
s}\bigr)_{n\in\bZ}$ is a Riesz basis of $L_2(0,1)$; moreover, there exists $m=m(h,r)>0$ that gives a lower bound of $\sE_{\bla}$ for every $\bla\in\sL(h,r)$~\cite{Hriesz}. Therefore,
\[
      ((I+\sH_+) f,f) = \mathrm{V.p.}\sum_{n=-\infty}^\infty \al_n |(f, \re^{2\la_n\ri s})|^2
            \ge h' \mathrm{V.p.}\sum_{n=-\infty}^\infty |(f, \re^{2\la_n\ri s})|^2
            \ge h'm\|f\|^2.
\]
The proof is complete.
\end{proof}

To study solvability of the Krein equation~\eqref{eq:inv.krein}, we shall regard it as a relation between the corresponding integral operators. To this end we recall several notions that will be used. By definition, the ideal $\fS_2= \fS_2(\bH)$ of Hilbert--Schmidt operators
consists of integral operators~$\sT$ in~$\bH$ whose kernels $T$ satisfy the inequality
\[
    \int_0^1\int_0^1  \tr \bigl(T(x,y)T^*(y,x)\bigr)\,dx\,dy <\infty;
\]
here $\tr$ stands for  the trace of a matrix or a trace-class operator. Notice that if $T=(T_{ij})_{i,j=1}^2$, then
\[
    \tr TT^* = \sum_{i,j=1}^2 |T_{ij}|^2,
\]
i.e., the operators in~$\fS_2$ have kernels whose entries are square integrable in~$\Omega= (0,1)\times(0,1)$. The linear set $\fS_2$ becomes a Hilbert space under the scalar product
\[
    \langle \sX,\sY\rangle_2 := \tr(\sX \sY^*)
        = \int_0^1\int_0^1 \tr \bigl(X(x,y)Y^*(y,x)\bigr)\,dx\,dy
\]

As noticed above, the matrix-valued function~$H$ satisfies the relation~$H^*(s)= H(-s)$ and $H(s)$ is unitarily equivalent to the diagonal matrix~$\diag\{h(s),h(-s)\}$. Therefore,
\begin{align*}
  \int_0^1 \tr H(x-y)H^*(x-y)\,dx
    &= \int_0^1 |h(x-y)|^2\,dx +  \int_0^1 |h(y-x)|^2\,dx \\
     &\le 2 \int_{-1}^{1} |h(s)|^2\,ds,
\end{align*}
so that the convolution operator~$\sH$ belongs to~$\fS_2$ and its norm satisfies $\|\sH\|^2_{\mathfrak S_2}=\langle \sH,\sH\rangle_2\le
4\|h\|^2_{L_2(0,2)}$.

Denote by $\fS_2^+$ the subspace of $\fS_2$ consisting of all
Hilbert--Schmidt operators with lower-triangular kernels. In other
words, $\sT \in \fS_2$ belongs to $\fS_2^+$ if the kernel $T$ of~$\sT$
satisfies $T(x,y) = 0$ for $0 \le x < y \le 1$. For an arbitrary
$\sT\in\fS_2$ with kernel $T$ the cut-off $T^+$ of $T$ given by
\[
    T^+(x,y) = \left\{ \begin{array}{ll}
        T(x,y)& \quad \mbox {for $x \ge y$}\\
              0 & \quad \mbox {for $x < y$} \end{array}
        \right.
\]
generates an operator $\sT^+\in\fS_2^+$, and the corresponding mapping
$\fP^+: \sT \mapsto \sT^+$  turns out to be an orthoprojector in
$\fS_2$ onto $\fS_2^+$, i.e. $(\fP^+)^2 = \fP^+$ and
 $\langle \fP^+ \sX, \sY\rangle_2 = \langle \sX,\fP^+ \sY\rangle_2$
for any $\sX,\sY\in \fS_2$; see details in~\cite[Ch.~I.10]{GK2}.

With these notations, the Krein equation~\eqref{eq:inv.krein} for the kernels $R$ and $H$ can be recast as the relation between the corresponding Hilbert--Schmidt integral operators
\begin{equation}\label{eq:KGLM}
    \sR + \fP^+ \sH + \fP^+(\sR\sH) = 0
\end{equation}
or
\[
    (\fI + \fP^+_{\sH}) \sR = - \fP^+ \sH,
\]
where $\fP^+_{\sX}$ is the linear operator in $\fS_2$ defined by
$\fP^+_{\sX} \sY = \fP^+ (\sY\sX)$ and $\fI$ is the identity operator
in $\fS_2$. Therefore solvability of the Krein equation and
continuity of its solutions on~$\sH$ is strongly connected with
the properties of the operator~$\fP^+_{\sH}$.


\begin{lemma}\label{lem:Pcnts}
For every $\sX\in\sB(\cH)$, the operator $\fP^+_\sX$ is bounded in~$\fS_2$. Moreover, for every $\sH$ from the set
\[
    \mathfrak{H} := \{ \sH = \sH(\bla,\bal) \mid (\bla,\bal) \in \sL(h,r)\times \sA(h',r')\} \subset \fS_2
\]
the operator $\fI + \fP^+_{\sH}$ is invertible in $\sB(\fS_2^+)$ and the inverse $(\fI + \fP^+_{\sH})$ depends analytically and Lipschitz continuously on $\sH\in\mathfrak{H}$ in the topology of~$\fS_2$.
\end{lemma}

\begin{proof}
Boundedness of $\fP^+_\sX$ is a straightforward consequence of the inequality
\[
    \|\fP^+_\sX \sY\|_{\fS_2} \le \|\sY\sX\|_{\fS_2}
        \le \|\sX\|_{\sB(\bH)} \|\sY\|_{\fS_2},
\]
cf.~\cite[Ch.~3]{GK1}. Assume next that $\sI + \sX \ge \eps \sI$ in~$\bH$; then for
$\sY\in\fS_2^+$ we find that
\[
    \langle (\fI + \fP^+_\sX) \sY, \sY\rangle_2
        = \langle \sY,\sY\rangle_2 + \langle \sY\sX,\sY \rangle_2
        = \tr \bigl(\sY(\sI+\sX)\sY^*\bigr).
\]
Since $\sY(\sI+\sX)\sY^* \ge \eps \sY\sY^*$ and the trace is a monotone functional, we
get
\[
    \langle (\fI+\fP^+_\sX) \sY,\sY \rangle_2 \ge
        \eps \langle \sY, \sY \rangle_2,
\]
i.e., $\fI + \fP^+_\sX \ge \eps \fI$ in $\fS_2^+$.

Applying now Lemma~\ref{lem:5.Iholds}, we conclude that for every $\sH\in \mathfrak{H}$ it holds $\fI + \fP^+_\sH \ge \eps \fI$ with $\eps$ of that lemma depending only on $h$, $h'$, $r$, and $r'$; therefore, $\fI + \fP^+_\sH $ is boundedly invertible in~$\sB(\fS_2^+)$ and
\[
    \bigl\|(\fI + \fP^+_\sH)^{-1} \bigr\| \le \eps^{-1}.
\]
Since $\fP^+_{\sH}$ depends linearly on~$\sH$, it follows that the mapping $\sH \mapsto (\fI + \fP^+_\sH)^{-1}$ from $\fS_2$ into $\sB(\fS_2^+)$ is analytic and Lipschitz continuous on the set~$\mathfrak H$. The proof is complete.
\end{proof}


\begin{corollary}\label{cor:inv.RonH}
For every $\sH\in\fH$, the Krein equation~\eqref{eq:inv.krein} has a unique solution
\[
    \sR:= - (\fI +\fP_\sH^+)^{-1} \fP^+ \sH \in \fS_2^+;
\]
moreover, $\sR$ depends analytically and Lipschitz continuously in~$\fS_2^+$ on $\sH \in \fH \subset \fS_2$.
\end{corollary}

It follows that the kernel~$R(x,t)$ of $\sR$ is square integrable in the domain~$\Omega$ and depends analytically and Lipschitz continuously in $L_2(\Omega)\otimes\bC^4$ on $\sH$. However, we need to know that $R(x,0)$ is well defined and belongs to $L_2(0,1)\otimes \bC^4$.

To this end we use the Krein equation to show that
\[
    R(x,t) = - H(x-t) - \int_0^1 R(x,s)H(s-t)\,ds
\]
as a function of~$x$ depends continuously in $L_2(\Omega) \times \bC^4$ on $t\in [0,1]$. Since the shift $f(\cdot)\mapsto f(\cdot-t)$ is a continuous operation in~$L_2(\bR)$, $H(\cdot-t)$ enjoys the required property. Denote by $R_{ij}$ and $H_{ij}$, $i,j =1,2$, the entries of the matrix-valued functions $R$ and $H$ respectively. They belong to $L_2(\Omega)$ since the corresponding integral operators $\sR$ and $\sH$ are in~$\fS_2$; therefore,
\begin{equation}\label{eq:inv.s2}
\begin{aligned}
    \int_0^1 \Bigl|\int_0^1 R_{ij}(x,s)& H_{jk}(s-t)\,ds\Bigr|^2\,dx\\
        &\le \int_0^1 dx \int_0^1 |R_{ij}(x,s)|^2\,ds
                        \int_0^1 |H_{jk}(s-t)|^2\,ds\\
        &\le \int_{-1}^1 |H_{jk}(s)|^2\,ds
            \int_0^1 \int_0^1 |R_{ij}(x,s)|\,ds\,dx < \infty.
\end{aligned}
\end{equation}
Thus the function
\begin{equation}\label{eq:inv.RH}
    \int_0^1 R_{ij}(x,s) H_{jk}(s-t)\,ds
\end{equation}
of the variable~$x\in[0,1]$ belongs to $L_2(0,1)$; moreover, continuity of the shifts $H_{jk}(\cdot-t)$ and estimate~\eqref{eq:inv.s2} show that function~\eqref{eq:inv.RH} depends continuously in $L_2(0,1)$ on $t\in[0,1]$. Since the entries of the matrix~$\int_0^1 R(x,s)H(s-t)\,ds$ are sums of functions~\eqref{eq:inv.RH} for corresponding indices~$i$, $j$ and $k$, we conclude that $R(\cdot,t)$ indeed depends continuously in $L_2(0,1)\otimes \bC^4$ on $t\in[0,1]$. In particular, $R(\,\cdot\,,0)$ is a well-defined function in $L_2(0,1)\otimes \bC^4$.

\begin{proof}[Proof of Theorem~\ref{thm:pre.norm}]
The above arguments show that it is legitimate to take $t=0 $ in the Krein equation; thus we get that
\[
    R(x,0) = - H(x) - \int_0^1 R(x,s)H(s)\,ds.
\]
The right-hand side is a bilinear expression in $H$ and $R$. In view of the analytic dependence of $\sR$ onto $\sH$ stated in Corollary~\ref{cor:inv.RonH} and estimates~\eqref{eq:inv.s2}, this yields analyticity and Lipschitz continuity of $R(x,0)$ on $H\in L_2(-1,1)\otimes \bC^4$. By~\eqref{eq:inv.QR} the function~$Q$ enjoys the same properties, and the proof is complete.
\end{proof}


\section{Reconstruction from two spectra}\label{sec:two}

We recall that the norming constants $\al_n$ for the Dirac operator~$\sD_1(Q)$ can be determined from the spectra~$(\la_n)$ and $(\mu_n)$ of $\sD_1(Q)$ and $\sD_2(Q)$ by the formula
\[
    \al_n = \frac1{\dot{S}(\la_n)C(\la_n)},
\]
where the entire functions~$S$ and $C$ are given by the canonical products~\eqref{eq:pre.PhiPsi} over~$\la_n$ and $\mu_n$ respectively. This induces a mapping $\bnu \mapsto \bal$ from the spectral data~$\bnu:=\bigl((\la_n),(\mu_n)\bigr)\in\sN$ into the norming constants~$\bal:=(\al_n)\in\sA$. In this section, we shall establish Theorem~\ref{thm:pre.main} by proving the following result.

\begin{theorem}\label{thm:two.norm}
For every $h\in(0,\pi/2)$ and $r>0$, the mapping
\begin{equation}\label{eq:two.nu-to-al}
   \sN(h,r) \ni \bnu \mapsto \bal \in \sA
\end{equation}
is analytic and Lipschitz continuous; moreover, there exist positive constants $h'$ and $r'$ such that the
range of this mapping belongs to~$\sA(h',r')$.
\end{theorem}

By definition, $\sA$ consists of elements of the commutative unital Banach algebra~$A$ introduced in Appendix~\ref{app:BA}. We observe that the metrics on~$\sA$ agrees with the norm of~$A$, and thus the results of Appendix~\ref{app:BA} yield the following statement.

\begin{proposition}\label{pro:two.inv}
For every positive $h$ and $r$, the set~$\sA(h,r)$ consists of invertible elements of~$A$. Moreover, the mapping~$\bal\mapsto\bal^{-1}$ is analytic and Lipschitz continuous in~$A$ on $\sA(h,r)$, and its range lies in $\sA\bigl((1+r)^{-1},rh^{-1}\bigr)$.
\end{proposition}

In view of Proposition~\ref{pro:two.inv}, it suffices to prove Theorem~\ref{thm:two.norm} with $\bal$ replaced by~$\bal^{-1}$. The elements of the sequence~$\bal^{-1}$ are $\al_n^{-1}=\dot{S}(\la_n)C(\la_n)$. We shall show that the sequences
\[
    \bga:=\bigl((-1)^n\dot{S}(\la_n)\bigr)_{n\in\bZ},
        \qquad
    \bde:=\bigl((-1)^n     C (\la_n)\bigr)_{n\in\bZ}
\]
form elements of~$\sA$. Thus Theorem~\ref{thm:two.norm} will be proved if we show that the mappings
\begin{equation}\label{eq:two.gade}
    \sN(h,r) \ni \bnu \mapsto \bga\in \sA, \qquad
    \sN(h,r) \ni \bnu \mapsto \bde\in \sA
\end{equation}
enjoy the  properties required for the mapping~\eqref{eq:two.nu-to-al}.

To begin with, integral representation~\eqref{eq:pre.sR} of the solution~$\bs(\cdot,\la)$ yields the formulae
\begin{align}\label{eq:two.S}
    S(\la) &= \sin\la + \int_0^1 r_1(s) \re^{\ri\la(1-2s)}\,ds,\\
    C(\la) &= \cos\la + \int_0^1 r_2(s) \re^{\ri\la(1-2s)}\,ds \label{eq:two.C}
\end{align}
for the functions~$S$ and $C$. Here $r_1$ and $r_2$ are some (uniquely defined) functions in~$L_2(0,1)$ formed from the entries of the matrix-valued function~$R(1,1-x)$. Therefore both $\dot{S}$ and $C$ can be recast in the form
\[
    \cos\la + \int_0^1 g(s) \re^{\ri\la(1-2s)}\,ds
\]
with $g(s)=\ri(1-2s)r_1(s)$ for $\dot{S}$ and $g(s) = r_2(s)$ for~$C$. The sequences $\bga$ and $\bde$ have therefore similar structures; namely, their $n$-th element equals
\begin{equation}\label{eq:two.ga-de-n}
    \cos\rho_{2n} + (-1)^n \int_0^1 g(s) \re^{\ri\la_n(1-2s)}\,ds
\end{equation}
for respective~$g$; here, as usual, $\rho_{2n}:= \la_n - \pi n$.

Clearly, the mapping $(\rho_{2n})\mapsto (\cos\rho_{2n}-1)$ is analytic in~$\ell_2$. Its Lipschitz continuity follows from the inequality $|\cos x - \cos y| \le |x-y|$; also, the inequality $1-\cos x \le x^2/2$ yields the estimate $\|(\cos\rho_{2n}-1)\|\le \|(\rho_{2n})\|^2/2$.

Let $f_{\bla}$ denote the function in~$L_2(0,1)$ such that $\hat f_{\bla}(n) = \rho_{2n}$ and set
\begin{equation}\label{eq:two.h}
    h_n:=(-1)^n \int_0^1 g(s) \re^{\ri\la_n(1-2s)}\,ds.
\end{equation}
Then the function~$h$ given by the Fourier series~$\sum_{n\in\bZ}h_n\re^{2\pi\ri ns}$ coincides with the function~$\Psi(f_{\bla},g)$ of Lemma~\ref{lem:B.Psi}.  It follows from that lemma that the sequence $(h_n)_{n\in\bZ}$ depends analytically and boundedly Lipschitz continuously in~$\ell_2$ on~$f_{\bla}$ and~$g$.
We prove in the lemma below that the function $g$ depends in the same manner on $\bnu = (\bla,\bmu) \in \sN(h,r)$.

\begin{lemma}\label{lem:two.kk1}
The mappings
\[
    \sN(h,r) \ni (\bla,\bmu)  \mapsto r_j \in L_2(0,1),
    \qquad j=1,2,
\]
are analytic and Lipschitz continuous.
\end{lemma}

\begin{proof}
Since both mappings can be treated similarly, we only discuss
the first one. By definition, we have $S(\la_n)=0$, and thus the numbers $\la_n = \pi n +
\rho_{2n}$, $n\in\bZ$, are zeros of the entire function~$S$ of~\eqref{eq:two.S}. The required properties of the mapping $\bnu \mapsto r_1$ follow now from the results of Appendix~\ref{app:map}.
\end{proof}


The above reasoning justifies the inclusion $\bal^{-1}\in\sA$ as well as analyticity and Lipschitz continuity of the mappings of~\eqref{eq:two.gade}. It remains to prove that there exist positive $h'$ and $r'$ such that, for every~$\bnu\in\sN(h,r)$, the corresponding elements $\bga$ and $\bde$ belong to $\sA(h',r')$.

Existence of such an $r'$ follows from the uniform estimates of the $\ell_2$-norms of the sequences $(\cos\rho_{2n}-1)$ and $(h_n)$; see~\eqref{eq:two.ga-de-n} and \eqref{eq:two.h}. Indeed, in view of Lemma~\ref{lem:B.Psi} the function $h=\sum_{n\in\bZ}h_n\re^{2\pi\ri ns}$ remains in the bounded subset of~$L_2(0,1)$ when $f_{\bla}$ and $g$ vary over bounded subsets of~$L_2(0,1)$, and the latter is the case when $\bnu$ runs over~$\sN(h,r)$ by the definition of the functions $f_{\bla}$ and $g$ and Lemma~\ref{lem:two.kk1}.

Next, in view of formula~\eqref{eq:pre.PhiPsi} and the interlacing property of~$\la_n$ and $\mu_n$, the numbers $\ga_n=(-1)^n\dot{S}(\la_n)$ and $\delta_n=(-1)^nC(\la_n)$ are all of the same sign and thus are all positive in view of the asymptotic relation~\eqref{eq:two.ga-de-n}. The uniform positivity of $\ga_n$ and $\delta_n$ (and thus existence of a positive $h'$ such that $1/\al_n=\ga_n\delta_n\ge h'$), immediately follows from the lemma below. A similar statement is established in~\cite{Hriesz}; we give a complete proof below for the sake of completeness.

\begin{lemma}\label{lem:two.PhiPsi}
For every $h\in (0,\pi/2)$ and $r>0$ we have
\[
    \sup_{(\bla,\bmu)} \sup_{n\in\bN}\,
           \log|\dot{S}(\la_n)|<\infty,
    \qquad
    \sup_{(\bla,\bmu)} \sup_{n\in\bN}\,
            \log|C(\la_n)| <\infty,
\]
where $S$ and $C$ are constructed via~\eqref{eq:pre.PhiPsi} from
sequences $\bla$ and $\bmu$, and the
suprema are taken over $(\bla,\bmu)\in \sN(h,r)$.
\end{lemma}

\begin{proof}
We assume first that $n\ne0$. By~\eqref{eq:pre.PhiPsi}, we have%
\begin{footnote}
{In what follows, all summations and multiplications over the index set
$\bZ$ will be taken in the principal value sense and the
symbol~$\mathrm{V.p.}$ will be omitted.}
\end{footnote}
\[
    \dot{S}(\la_n) = \frac{\la_0-\la_n}{\pi n}\prod_{k\in\bZ,\ k\ne 0,n}\frac{\la_k - \la_n}{\pi k}.
\]
Dividing both sides by
\[
    \cos \pi n = \frac{d \sin z}{dz}\Bigr|_{z=\pi n}
                = - \prod_{k\in\bZ,\ k\ne 0,n}
                    \frac{\pi k -\pi n}{\pi k},
\]
we conclude that

\[
    |\dot{S}(\la_n)| =
            \prod_{k\in\bZ, \ k\ne n} \frac{\la_k - \la_n}{\pi (k-n)}
\]
and observe that this formula also holds for $n=0$. Set (recall that $\rho_{2k}:=\la_k-\pi k$)
\[
    a_{k,n}:= \frac{\la_{k} - \la_{n}}{\pi(k-n)} -1
            = \frac{\rho_{2k}-\rho_{2n}}{\pi (k-n)}
\]
if $k\ne n$ and $a_{n,n}:=0$; then
    \(
    |\dot{S}(\la_n)| = \prod_{k\in\bZ}
        (1+ a_{k,n}).
    \)
Since the sequence $(\la_n)$ is $2h$-separated for every
$(\bla,\bmu)\in\sN(h,r)$, we have
 \(
    1+a_{k,n} \ge {2h}/\pi
 \)
for all integer $k$ and~$n$. Therefore, with
\[
    K:= \max_{x\ge -1 + 2h/\pi} \Bigl|\frac{\log(1+x)-x}{x^2}\Bigr|
    < \infty,
\]
we get the estimate
\begin{equation}\label{eq:two.log}
    \Bigl|\log \prod_{k\in\bZ}(1+a_{k,n})\Bigr|
            \le \Bigl|\sum_{k\in\bZ} a_{k,n}\Bigr|
            + K \sum_{k\in\bZ} a^2_{k,n},
\end{equation}
provided the two series converge.

Clearly,
\[
    \sum_{k\ne n} \frac{1}{k-n}  =0,
\]
and thus
\[
    \Bigl|\sum_{k\in\bZ} a_{k,n}\Bigr|
        = \Bigl|\frac1\pi\sum_{k\ne n} \frac{\rho_{2k}}{k-n}\Bigr|
        \le \frac{r}{\sqrt3}
\]
by the Cauchy--Bunyakovski--Schwarz inequality (recall that
$\sum_{k\in\bZ}\rho^2_{2k} \le r^2$ by the definition of the set
$\sN(h,r)$ and $\sum_{k\ne n}(k-n)^{-2}=\pi^2/3$). Next, the
inequality
\[
    a_{k,n}^2 \le  \frac{2 \rho_{2k}^2}{(k-n)^2}
                   + \frac{2 \rho_{2n}^2}{(k-n)^2}
\]
for $k\ne n$ yields
\[
    \sum_{k\in\bZ} a_{k,n}^2
        \le 4r^2 \sum_{k\ne n}\frac{1}{(n-k)^2}
        = \frac{4\pi^2r^2}3.
\]
It follows from~\eqref{eq:two.log} that
\[
    \Bigl|\log \prod_{k\in\bZ}(1+a_{k,n})\Bigr|
        \le (\sqrt3 r + 4K \pi^2r^2)/3,
\]
where the constant~$K$ only depends on $h$.

Similarly, we find that
\[
    |{C(\la_n)}|
        = \Bigl|\prod_{k=1}^\infty
        \frac{\mu_k-\la_n}{\pi(k+\tfrac12)}\Bigr|
        = \prod_{k\in\bZ}
        \frac{\mu_k -\la_{n}}{\pi(k+\tfrac12) - \pi n}
\]
and then mimic the above reasoning to establish the other uniform bound.
The lemma is proved.
\end{proof}

Combination of the above results constitutes a complete proof of Theorem~\ref{thm:two.norm} and together with Theorem~\ref{thm:pre.norm} establishes Theorem~\ref{thm:pre.main}.

\medskip
\textbf{Acknowledgements.} {The author thanks A.~A.~Shkalikov and Ya.~V.~Mykytyuk for stimulating
discussions. The research was partially supported by the Alexander von
Humboldt Foundation and was partially carried out during the visit to the
Institute for Applied Mathematics of Bonn University, whose hospitality
is warmly acknowledged.}

\appendix

\section{Analyticity of some related mappings}\label{app:map}

Here we give a brief account on the results of~\cite{HMzero} as well as some extensions that are used to prove Lemma~\ref{lem:two.kk1}. It was shown in~\cite{HMzero} that for every $f\in L_2(0,1)$ there exists a unique function~$g\in L_2(0,1)$ such that
all zeros (counting multiplicities) of the entire function
\begin{equation}\label{eq:A1.g}
    G_g(z):=    \sin z + \int_0^1 g(t) \re^{\ri z(1-2t)}\,dt
\end{equation}
are given  by the numbers $\pi n + \hat
f(n)$, $n\in\bZ$. Such pairs of $f$ and $g$ in fact satisfy the
relation
\begin{equation}\label{eq:2.H}
    H(f,g) := s(f) + g +
        \sum_{k=1}^\infty \frac{(M^k g)\ast f^{<k>}}{k!} =0;
\end{equation}
here
\[
    s(f):= \sum_{k=0}^\infty\frac{(-1)^k f^{<2k+1>}}{(2k+1)!},
\]
$f^{<k>}$ is the $k$-fold convolution of~$f$ with itself, and $M$ is
the operator of multiplication by~$\ri (1-2x)$. The function $H$ is
analytic from $L_2(0,1)\times L_2(0,1)$ into $L_2(0,1)$, and its
partial derivatives $\partial_f H(f,g)$ and $\partial_g H(f,g)$ are
given by
\begin{align}\label{eq:5.partialf}
    \partial_f H(f,g)(h_1) &=
         \Bigl( c(f) + \sum_{k=1}^\infty \frac{(M^k g)\ast f^{<k-1>}}{(k-1)!}
         \Bigr) \ast h_1,\\
    \partial_g H(f,g)(h_2) &=
        h_2 + \sum_{k=1}^\infty \frac{(M^k h_2)\ast f^{<k>}}{k!}
        \label{eq:5.partialg}
\end{align}
with
\[
    c(f) := \sum_{k=0}^\infty \frac{(-1)^k f^{<2k>}}{(2k)!}.
\]
Using the implicit function theorem, it was shown that the induced mapping $f \mapsto g$ is analytic. In order to establish its Lipschitz continuity, we shall study the above partial derivatives in more detail.

Namely, we assume that $f\in L_2(0,1)$ is such that the
corresponding sequence $\bla=(\la_n)_{n\in\bN}$ with $\la_n := \pi n + \hat f(n)$ belongs to $\sL(h,r)$. Set
$S_{\bla}$ to be the canonical product of~\eqref{eq:pre.PhiPsi}; then $S_{\bla}$ can also be represented as~\eqref{eq:A1.g}. Direct calculations show that the $n$-th Fourier coefficient of the function
of~\eqref{eq:5.partialf} is equal to
\[
    (-1)^n \hat h_1(n) \Bigl[ \cos \la_n
            + \int_0^1 \ri (1-2t)g(t) \re^{\ri \la_n(1-2t)}\,dt \Bigr]
            = (-1)^n \hat h_1(n) \dot{S}_{\bla}(\la_n).
\]
By Lemma~\ref{lem:two.PhiPsi} there are positive numbers $K_1$ and $K_2$ such
that
\[
    K_1 \le |\dot{S}_{\bla}(\la_n)| \le K_2
\]
for all $\bla\in \sL(h,r)$ and all $n\in\bZ$. Therefore the partial derivative $\partial_f H(f,g)$ is a bounded and boundedly invertible operator in $L_2(0,1)$; moreover, for every fixed $h>0$ and $r>0$, the norms of $\partial_f H(f,g)$ and their inverses are uniformly bounded for $f\in L_2(0,1)$ generating the sequences $\bla$ in the set $\sL(h,r)$.

Similarly, the $n$-th Fourier coefficient of the function
of~\eqref{eq:5.partialg} is equal to
\[
    (-1)^n \int_0^1 h_2(t) \re^{\ri \la_n(1-2t)}\,dt.
\]
By the results of~\cite{Hriesz}, there exist positive $M$ and $m$ such
that, for all $\bla\in \sL(h,r)$, the sequences
 $(\re^{\ri\la_n(1-2x)})_{n\in\bZ}$
form Riesz bases of $L_2(0,1)$ of upper bound~$M$ and lower
bound~$m$. Therefore the
operator~$H_g:=\partial_g H(f,g)$,
\[
   H_g :\,h_2 \mapsto
        \sum_{n\in\bZ} (-1)^n(h_2, \re^{\ri\la_n(1-2x)})
            \,\re^{2\pi n\ri x},
\]
is bounded and boundedly invertible in~$L_2(0,1)$, with $\|H_g\|\le
M^{1/2}$ and $\|H_g^{-1}\|\le m^{-1/2}$.

We now use the implicit mapping theorem to conclude that the mapping
$f\mapsto g$ is analytic in $L_2(0,1)$. The uniform bounds on the inverses of the partial
derivatives $\partial_f H(f,g)$ and $\partial_g H(f,g)$ established
above imply that, for every $h>0$ and $r>0$, this mapping is Lipschitz continuous on the set of functions~$f\in L_2(0,1)$ generating the sequences $\bla\in\sL(h,r)$.

\section{Some auxiliary results}\label{app:aux}

We recall that the convolution $f\ast g$ of two functions in~$L_2(0,1)$ is a function in $L_2(0,1)$ given by
\[
    (f\ast g) (x) := \int_0^1 f(x-t) g(t)\,dt,
\]
where $f$ is extended to $(-1,0)$ as a periodic function with period~$1$.
The (discrete) Fourier transform $\hat f$ of $f\in L_2(0,1)$ is a function over $\bZ$ given by
\[
    \hat f(n) := \int_0^1 f(t) \re^{-2\pi n \ri t}\,dt.
\]
It is well known that the Fourier transform is a unitary mapping from $L_2(0,1)$ to $\ell_2(\bZ)$ and that
$\widehat{f\ast g}(n) = \hat f(n) \hat g(n)$; as a result, we have the inequality
\[
    \|f\ast g\| \le \|f\|\|g\|
\]
for all $f,g\in L_2(0,1)$.

\begin{lemma}\label{lem:B.gf}
For a function $f\in L_2(0,1)$, set
\[
    g(f)(x) := \mathrm{V.p.}\sum_{n=-\infty}^\infty
        [\re^{2\hat f(n) \ri x} - 1]\re^{2\pi n \ri x}.
\]
Then the series determines a function in~$L_2(0,1)$, and the mapping
\[
    L_2(0,1) \ni f \mapsto g_f \in L_2(0,1)
\]
is analytic and locally Lipschitz continuous on bounded subsets.
\end{lemma}

\begin{proof}
We start with observing that the series
 \(
    \sum_{n\in\bZ} \hat f^{\,k}(n) \re^{2\pi n\ri s}
 \)
is the Fourier series for the function~$f^{\langle k\rangle}$, the $k$-fold
convolution of~$f$ with itself, and that $\|f^{\langle k\rangle}\| \le \|f\|^k$.  Developing $\re^{\hat f(n) \ri s}$ into the Taylor series, we find that
\begin{align*}
    g(f)  &= \mathrm{V.p.}\sum_{n=-\infty}^\infty
        \Bigl[\sum_{k=1}^\infty \frac{\hat f^k(n) (2\ri s)^k}{k!} \Bigr]
            \re^{2\pi n\ri s}\\
        &= \sum_{k=1}^\infty \frac{(2\ri s)^k}{k!}
            \mathrm{V.p.}\sum_{n=-\infty}^\infty\hat f^{k}(n)
            \re^{2\pi n\ri s}\\
        &= \sum_{k=1}^\infty \frac{(2\ri s)^k}{k!} f^{\langle k \rangle}.
\end{align*}
The change of the summation order in the second equality above is
justified by the fact that, for $k>1$, the summands in the double
series are dominated by $C^k\hat f^2(n)/k!$ with
$C:=2\max_{n\in\bZ}\{|\hat f(n)|\}+1$. Therefore the double series over the
index set $\{(n,k) \mid n\in\bZ, k>1\}$ converges absolutely and the
Fubini theorem applies. This formula represents $g(f)$ as an absolutely convergent series (which is a Taylor series expansion of~$g(f)$ in the variable~$f$) and thus proves
the analyticity in~$L_2(0,1)$ of the mapping~$f\mapsto g(f)$.

Lipschitz continuity of that mapping on bounded sets follows from the estimate
\begin{align*}
    \|g(f_1) - g(f_2)\|
            &= \Bigl\|\sum_{k=1}^\infty \frac{(2\ri s)^k}{k!}
            [f_1^{\langle k \rangle} -  f_2^{\langle k \rangle}]\Bigr\|\\
     &\le \sum_{k=1}^\infty \frac{2^k}{(k-1)!}  \|f_1-f_2\|
                    \bigl( \|f_1\| + \|f_2\| \bigr)^{k-1}
                    \le \exp\{4r\} \|f_1- f_2\|,
\end{align*}
which is valid as soon as the $L_2$-norms of $f_1$ and $f_2$ are not greater than~$r$. The proof is complete.
\end{proof}

\begin{lemma}\label{lem:B.Phi}
For $f$ and $g$ in~$L_2(0,1)$, set
\[
    \Phi(f,g):= \mathrm{V.p.}\sum_{n=-\infty}^\infty
        \hat g(n) \exp\{2[\pi n +\hat f(n)] \ri s\}.
\]
Then the function $\Phi(f,g)$ belongs to~$L_2(0,1)$ and the mapping
\[
    \Phi\,:\, L_2(0,1)\times L_2(0,1) \to L_2(0,1)
\]
is analytic and Lipschitz continuous on bounded subsets.
\end{lemma}

\begin{proof}
Transformations similar to those used in the proof of the above lemma show that
\[
    \Phi(f,g) = \sum_{k=1}^\infty \frac{(2\ri s)^k}{k!} \, [f^{\langle k \rangle} \ast g].
\]
The mapping~$\Phi$ is linear (and thus analytic) in $g$, and its analyticity in~$f$
as well as Lipschitz continuity on bounded subsets is established in the same manner as for the mapping $g(f)$ of Lemma~\ref{lem:B.gf}.
\end{proof}

\begin{lemma}\label{lem:B.Psi}
For $f$ and $g$ in~$L_2(0,1)$, set
\[
    \Psi(f,g):= \mathrm{V.p.}\sum_{n\in\bZ}
        (-1)^n \int_0^1 g(t) \exp\{[\pi n+ \hat f(n)]\ri(1-2t)\}\,dt\, \re^{2\pi\ri nx};
\]
Then the function $\Psi(f,g)$ belongs to~$L_2(0,1)$ and the mapping
\[
    \Psi\,:\, L_2(0,1)\times L_2(0,1) \to L_2(0,1)
\]
is analytic and Lipschitz continuous on bounded subsets.
\end{lemma}

\begin{proof}
The coefficient of $\re^{2\pi\ri nx}$ in the above series for~$\Psi$
can be written as
\begin{equation}\label{eq:C.c2n}
    \int_0^1 g(t) \exp\{ \ri(1-2t)\hat f(n)\} \re^{-2\pi\ri nt}\,dt
\end{equation}
and gives the $n$-th Fourier coefficient of the function
 \(
    h:=\sum_{k=0}^\infty {h_k}/{k!},
 \)
with $h_0:=g$, $h_k:= f^{<k>} \ast M^k g$ for $k\ge1$, and $M$
being the operator of multiplication by the function~$\ri(1-2t)$.
In other words, we have $\Psi(f,g)=h$. Since $\|f\ast g\|\le \|f\|\|g\|$ for every $f$ and $g$ in~$L_2(0,1)$, the functions $h_k$ belong to $L_2(0,1)$ and their norms there obey the estimate
\[
    \|h_k\|  \le  \|f\|^k \|M\|^k\|g\|,
\]
with $\|M\|$ denoting the norm of the operator $M$. Thus the series for~$h$ converges absolutely and, since every $h_k$ is a multi-linear function of~$f$ and $g$, the mapping $\Psi$ is analytic. Its Lipschitz continuity on bounded subsets is established in the usual manner, and the proof is complete.
\end{proof}

\section{Banach algebras}\label{app:BA}

The space $\ell_2=\ell_2(\bZ)$ is a commutative Banach algebra under the pointwise multiplication $(x_n)\cdot (y_n) = (x_n\cdot y_n)$. Its unital extension~$A$ consists of elements of~$\ell_\infty$ of the form~$a\mathbf1 + \bx$ with $a\in\bC$, the unity~$\mathbf 1\in\ell_\infty$ having all its elements equal to~$1$, and $\bx=(x_n)\in\ell_2$. The norm in~$A$ is given by
\[
    \|a\mathbf{1} + \bx\|_A = |a| + \|\bx\|,
\]
and $a\mathbf{1} +\bx$ is invertible in~$A$ if and only if $a\ne0$ and $a+x_n\ne0$ for all~$n\in\bZ$; in this case the inverse is equal to $a^{-1}\mathbf{1} + \by$, where $\by = (y_n)$ with $y_n:=-x_n/a(a+x_n)$.
Since under the above assumptions we have
$
    \inf_n |a + x_n| > 0
$, we see that $\by$ indeed belongs to~$\ell_2$; moreover,
\[
    \|(a\mathbf{1} +\bx)^{-1}\|_A \le
            |a|^{-1}\bigl(1 + \|\bx\|/\inf_n|a+x_n|\bigr).
\]
The mapping $\hat\bx \mapsto \hat\bx^{-1}$ is analytic on the open set of all invertible elements of~$A$; in addition, it is Lipschitz continuous on the sets
\[
    \sS_\eps:= \bigl\{a\mathbf{1} +\bx \mid
                    |a|\ge\eps, \ \inf_n|a+x_n|\ge\eps\bigr\}.
\]



\begin{thebibliography}{99}

\bibitem{AKNS} M.~J.~Ablowitz, D.~J.~Kaup, A.~C.~Newell, and
              H.~Segur,
    Nonlinear-evolution equations of physical significance,
    \emph{Phys. Rev. Lett.} {\bf 31} (1973), 125--127.

\bibitem{AKNS2} M.~J.~Ablowitz, D.~J.~Kaup, A.~C.~Newell, and H.~Segur,
    The inverse scattering transform---Fourier analysis for nonlinear problems,
    \emph{Stud. Appl. Math.}
    \textbf{53} (1974), 249--315.


\bibitem{AHMdir} S.~Albeverio, R.~Hryniv, and Ya.~Mykytyuk,
    Inverse spectral problems for Dirac operators with summable
    potentials,
    \emph{Russian J. Math. Phys.}
    \textbf{12} (2005), no.~4, 406--423; arXiv preprint
    math.SP/0701158.

\bibitem{Am} B.~A.~Ambartsumyan,
    \"Uber eine Frage der Eigenwerttheorie,
    \emph{Zeitschr. f\"ur Physik} \textbf{53}(1929), 690--695.

\bibitem{Amour} L.~Amour,
    Inverse spectral theory for the AKNS system with separated boundary conditions,
    \emph{Inverse Problems} \textbf{9} (1993), no.~5, 507--523.

\bibitem{An} L.~Andersson,
    Inverse eigenvalue problems for a Sturm--Liouville equation in impedance
    form,
    \emph{Inverse Probl.} \textbf{4} (1988), 929--971.


\bibitem{CG} S.~Clark and F.~Gesztesy,
    Weyl-Titchmarsh $M$-function asymptotics, local uniqueness results,
    trace formulas, and Borg-type theorems for Dirac operators,
    \emph{Trans. Amer. Math. Soc.} {\bf 354} (2002), no.~9,
    3475--3534.

\bibitem{CM} C.~F.~Coleman and J.~R.~McLaughlin,
    Solution of the inverse spectral problem for an impedance with
    integrable derivative,
    I, \emph{Comm. Pure Appl. Math.} \textbf{46} (1993), 145--184;
    II, \emph{Comm. Pure Appl. Math.} \textbf{46} (1993), 185--212.

\bibitem{GasL} M.~G.~Gasymov and B.~M.~Levitan,
    Determination of a differential operator from two spectra,
    \emph{Uspekhi Matem. Nauk} \textbf{19} (19964), no.~2, 3--63.

\bibitem{GasLscat} M.~G.~Gasymov\ and\ B. M. Levitan,
    Determination of the {D}irac system from the scattering phase,
    \emph{Dokl. Akad. Nauk SSSR} {\bf 167} (1966), 1219--1222.

\bibitem{GasLsp} M.~G.~Gasymov\ and\ B.~M.~Levitan,
    The inverse problem for the {D}irac system,
    \emph{Dokl. Akad. Nauk SSSR} {\bf 167} (1966), 967--970.

\bibitem{GL} I.~M.~Gelfand and B.~M.~Levitan,
    On determination of a differential equation by its spectral function,
    \emph{Izv. Akad. Nauk SSSR, Ser. Mat.} \textbf{15}(1951), no.~4,
    309--360 (in Russian).

\bibitem{GK1} I.~Gohberg and M.~Krein,
    \emph{Introduction to the Theory of Linear Non-selfadjoint Operators
    in Hilbert Space},
    Nauka Publ., Moscow, 1965 (in Russian); \emph{Engl. transl.:}
    Amer. Math. Soc. Transl. Math. Monographs, vol.~18,
    Amer. Math. Soc., Providence, RI, 1969.

\bibitem{GK2} I.~Gohberg and M.~Krein,
    \emph{Theory of Volterra Operators in Hilbert Space and its
    Applications},
    Nauka Publ., Moscow, 1967 (in Russian); \emph{Engl. transl.:}
    Amer. Math. Soc. Transl. Math. Monographs, vol.~24,
    Amer. Math. Soc., Providence, RI, 1970.

\bibitem{HJKS} D.~B.~Hinton, A.~K.~Jordan, M.~Klaus, and J.~K.~Shaw,
    Inverse scattering on the line for a {D}irac system,
    \emph{J. Math. Phys.} \textbf{32} (1991), no.~11, 3015--3030.

\bibitem{Ho} M.~Horv\'ath,
    On a theorem of Ambarzumyan,
    \emph{Proc. Roy. Soc. Edinb. A} \textbf{131} (2001), 899--907.

\bibitem{Hriesz} R.~O.~Hryniv,
    Uniformly bounded families of Riesz bases of exponentials, sines, and
    cosines,
    \emph{Mathem. Zametki} \textbf{87} (2010), no.~4, 542--553 (in Russian); \emph{Engl. transl.:} \emph{Mathem. Notes} \textbf{87} (2010), no.~4, 510--520.

\bibitem{Hstab-imp} R.~O.~Hryniv,
    Analyticity and uniform stability in the inverse spectral problem for impedance Sturm--Liouville operators,
    \emph{Carpat. Math. Publ.} \textbf{2} (2010), no.~1, 35--58.

\bibitem{Hstab} R.~O.~Hryniv,
    Analyticity and uniform stability of the inverse singular Sturm--Liouville spectral problem,
    \emph{submitted}; \texttt{arXiv:1101.5426v1 [math.SP]}.

\bibitem{HMinv}  R.~O.~Hryniv and Ya.~V.~Mykytyuk,
    Inverse spectral problems for Sturm--Liouville operators
    with singular potentials,
    \emph{Inverse Problems} \textbf{19} (2003), 665--684.

\bibitem{HMscale} R.~O.~Hryniv and Ya.~V.~Mykytyuk,
    Inverse spectral problems for Sturm--Liouville operators
    with singular potentials, IV. Potentials in the Sobolev space
    scale,
    \emph{Proc. Edinb. Math. Soc.}
    \textbf{49} (2006), no.~2, 309--329.


\bibitem{HMzero} R.~O.~Hryniv and Ya.~V.~Mykytyuk,
    On zeros of some entire functions,
    \emph{Trans. AMS}, \textbf{361} (2009), no.~4, 2207--2223.

\bibitem{Ki} M.~Kiss,
    An n-dimensional Ambarzumyan type theorem for Dirac operators,
    \emph{Inverse Problems} \textbf{20} (2004), 1593--1597.


\bibitem{LS} B.~M.~Levitan and I.~S.~Sargsjan,
    Sturm--Liouville and Dirac Operators,
    Nauka, Moscow, 1988 (in Russian);
    \emph{Engl. transl.:} Kluwer Acad., Dordrecht, 1991.

\bibitem{Mal99} M.~M.~Malamud,
    Questions of uniqueness in inverse problems for systems of differential equations on a finite
    interval,
    \emph{Tr. Mosk. Mat. Ob-va} {\bf 60} (1999), 199--258; translation in
    \emph{Trans. Moscow Math. Soc.} {\bf 1999}, 173--224.

\bibitem{Ma} V.~A.~Marchenko,
    \emph{Sturm-Liouville Operators and Their Applications},
    Naukova Dumka Publ., Kiev, 1977 (in Russian); \emph{Engl. transl.:}
    Birkh\"auser Verlag, Basel, 1986.

\bibitem{McL} J.~R.~McLaughlin,
    Stability theorems for two inverse problems,
    \emph{Inverse Probl.} \textbf{4}(1988), 529--540.

\bibitem{Rya3} T.~I.~Ryabushko,
    Estimation of the norm of the difference of two potentials of
    Sturm--Liouville boundary value problems,
    \emph{Teor. Funktsi\u\i\ Funktsional. Anal. i Prilozhen.}
    \textbf{39} (1983), 114--117 (in Russian).

\bibitem{SSstab} A.~M.~Savchuk and A.~A.~Shkalikov,
    Inverse problems for Strum--Liouville operators with potentials in Sobolev spaces: Uniform stability, \emph{Funktsion. Anal. Prilozhen.}
    \textbf{44} (2010), no.~4, 34--53 (in Russian); \emph{Engl. transl.:} \emph{Funct. Anal. Appl.} \textbf{44} (2010), no.~4, 270--285.

\bibitem{Th} B.~Thaller,
    {\it The Dirac Equation}, Springer, Berlin, 1992.


\bibitem{Yang} C.~F.~Yang and Z.~Y.~Huang,
    Inverse spectral problems for 2m-dimensional canonical Dirac operators,
    \emph{Inverse Problems} \textbf{23} (2007), 2565--2574.

\bibitem{Young} R.~Young,
    \emph{An Introduction to Nonharmonic Fourier Series},
    Academic Press, New York, (revised first edition), 2001.

\bibitem{ZS} V.~E.~Zakharov and A.~B.~Shabat,
    Exact theory of two-dimensional self-focusing and one-dimensional
    self-modulation of waves in nonlinear media,
    \emph{Soviet Physics JETP} {\bf 34} (1972), no.~1, 62--69.;
    translated from \emph{Zh. \`Eksper. Teoret. Fiz.}
    {\bf 61} (1971), no. 1, 118--134.

\end{thebibliography}
\end{document}